\documentclass{amsart}

\usepackage{amsmath}
\usepackage{amsfonts}
\usepackage{amsopn}
\usepackage{amssymb}

\newcommand{\remark}{\noindent \emph{Remark: }}
\newcommand{\R}{\ensuremath{\mathbb{R}}}

\newcommand{\Hyp}{\ensuremath{\mathbb{H}}}
\DeclareMathOperator{\Isom}{Isom}
\DeclareMathOperator{\length}{length}
\DeclareMathOperator{\Comm}{Comm}
\DeclareMathOperator{\diam}{diam}
\DeclareMathOperator{\vol}{vol}

\newtheorem{thm}{Theorem}
\newtheorem{lemma}[thm]{Lemma}

\title{Counting Hyperbolic Manifolds with Bounded Diameter}
\author{Robert Young}
\address{Department of Mathematics\\
University of Chicago\\
5734 S. University Avenue\\
Chicago, Illinois 60637\\
U.S.A.}
\keywords{diameter, hyperbolic manifolds}
\subjclass[2000]{Primary: 57N16; Secondary: 22E40}

\begin{document}
\bibliographystyle{plain}
\begin{abstract}
Let $\rho_n(V)$ be the number of complete hyperbolic manifolds of
dimension $n$ with volume less than $V$.  Burger, Gelander, Lubotzky,
and Moses\cite{BGLM} showed that when $n\ge 4$ there exist $a,b>0$
depending on the dimension such that $a V\log{V}\le \log{\rho_n(V)}\le
bV\log{V},$ for $V\gg 0$.  In this note, we use their methods to bound
the number of hyperbolic manifolds with diameter less than $d$ and
show that the number grows double-exponentially.
Additionally, this bound holds in dimension 3.
\end{abstract}
\maketitle
\section{Introduction}
In dimensions 4 and larger, the number of hyperbolic manifolds with a
bounded volume is finite.  Gelander\cite{Gelander04} has proven a similar finiteness result for manifolds locally isometric to symmetric spaces of dimension at least $4$, with the exceptions of $\Hyp^2\times \Hyp^2$ and $\text{SL}_3(\R)/\text{SO}_3(\R)$.  For $\Hyp^3$, Thurston\cite{Thurston} has
shown that this is not the case and that an infinite number of manifolds
of bounded volume can be constructed using Dehn surgery.  The
diameter of these manifolds, however, becomes large.  This note shows that the
number of hyperbolic manifolds of bounded diameter is finite and
provides upper and lower bounds on its growth.  We prove the following bounds.

\begin{thm}
Let $\tau_n(d)$ be the number of closed hyperbolic manifolds of
dimension $n$ with diameter less than or equal to $d$.  There exist
constants $a,b>0$ depending on the dimension such that for $d\gg 0$,
$$e^{e^{ad}}<\tau_n(d)<\begin{cases}e^{b d e^{(n-1) d}} & n\ge 4 \\
				      e^{b d e^{5d}} & n=3 \end{cases}.$$
\end{thm}
The author would like to thank Benson Farb and Shmuel Weinberger for
suggesting this question.

\section{Constructing the Lower Bound}
We achieve the lower bound by constructing exponentially
many(in the index of the cover) covers of a non-arithmetic hyperbolic
manifold with diameter growing logarithmically in the index.  We then
claim that only a fraction of these can be isomorphic.  This is a modification of the argument given in \cite{BGLM}.

We first note that the diameter of a cover can be estimated by its fundamental group as in the following lemma.

\begin{lemma}\label{diamlemma}
Let $M$ be a compact manifold, with $\Gamma=\pi_1(M)$, and
$\{\alpha_1,\dots,\alpha_n\}$ a set of generators.  Let $M'$ be a
finite cover of $M$ with $\pi_1(M')=\Gamma'$ and projection map
$p:M'\to M$.  There are constants $c_1, c_2$ depending on $M$ and the choice
of generating set such that if $g_1,\dots,g_r\in\Gamma$ are a set such
that $\cup_{i=1}^rg_i\Gamma'=\Gamma$ then 
$$\diam{M'}\le c_1+c_2\max\{w(g_i)\}_{i=1}^r,$$ 
where $w(g_i)$ is the word length of $g_i$ with respect to
$\{\alpha_1,\dots,\alpha_n\}$.
\end{lemma}

\proof{This follows from the observation that $\Gamma$ is quasi-isometric to
the universal cover of $M$.  Indeed, if we fix a basepoint $*\in M'$,
and loops $a_i:[0,1]\to M$ based at $*$ and representing the
$\alpha_i$, then one can choose $c_1=2\diam{M}$ and
$c_2=2\max\{\length(a_i)\}$.}

Gromov and Piatetski-Shapiro proved in \cite{GPS} that there exists a
non-arithmetic torsion-free cocompact lattice $\Gamma_n$ in
$\Isom{\Hyp^n}$ for $n\ge 2$.  Lubotzky\cite{Lubotzky96} showed that this group
contains a finite index subgroup $\Gamma\subset \Gamma_n$ such that there is a
surjection $f:\Gamma \to F_2$, where $F_2$ is the free group on two
generators.  We will use finite covers of the manifold $M=\Hyp^n/\Gamma$
to establish the lower bound.

Given an index $r$ subgroup $G$ of $F_2$, we can take the finite
sheeted cover $\widehat{M}$ of $M$ induced by the subgroup
$f^{-1}(G)\subset\Gamma$.  We can then use Lemma \ref{diamlemma} to estimate the
diameter of $\widehat{M}$.  Fix $\alpha_1, \alpha_2$
generators of $F_2$; we will construct  
many finite index subgroups of $F_2$ with sets of
relatively short coset representatives.  

We first give a correspondence between certain pairs of permutations of
$\{0,\dots,r-1\}$ and subgroups of $F_2$.  Given a pair of
permutations $(\sigma_1,\sigma_2)$, we can consider the action of
$F_2$ on $\{0,\dots,r-1\}$ where $\alpha_i$ acts by $\sigma_i$.  If
this action is transitive, the subset of $F_2$ that fixes $0$ is a
subgroup of index $r$, and the $i$th coset is the set of elements
taking $0$ to $i$.  In addition, if the pair $(\sigma'_1,\sigma'_2)$
determines the same subgroup as $(\sigma_1,\sigma_2)$, then there
exists a permutation leaving 0 fixed that conjugates
$\sigma'_1$ to $\sigma_1$ and $\sigma'_2$ to $\sigma_2$.  Conversely,
given a finite index subgroup of $F_2$, the action of $F_2$ on its
cosets by left multiplication gives such a pair of permutations,
defined up to relabeling the cosets that are not the original
subgroup.

Let $\sigma_1$ take $i$ to $i+1$ mod $r$.  Let $S$ be the set of
permutations $\sigma$ such that $\sigma(i)=2i$ for $i$ even, $0 < i <
r/2$, and $\sigma_2\in S$.  Then, if $d_k\dots d_1d_0$ is the binary expansion of $i < r/2$
and $d_k=1$, then
\begin{equation*}
\begin{split}
(\sigma_1^{2d_{0}}\circ \sigma_2 \circ \sigma_1^{2d_1}\circ &\dots\circ \sigma_2 \circ \sigma_1^{2d_{k-1}}\circ \sigma_2\circ \sigma_1^{2d_k})(0) \\
&=((\dots((0+2d_k)\cdot 2+2d_{k-1})\dots)\cdot 2+2 d_0) \\
&=2d_k\cdot2^k+2d_{k-1}\cdot2^{k-1}+\dots+2d_0 \\ 
&=2 i
\end{split}
\end{equation*}
This gives a coset representative for the $2i$th coset of length at
most $3(1+\log_2 i)$.  By composing this with $\sigma_1$, we obtain a
representative for the $(2i+1)$th coset.  Since
$|S|=(\lfloor{r/2}\rfloor+1)!$, we can use these to construct
$(\lfloor{r/2}\rfloor+1)!$ $r$-sheeted covers with diameter at most $D
\log{r}$, for some constant $D>0$.  Some of these, however, may be
isometric, but we can limit this effect by using the
non-arithmeticity of $\Gamma$.

Indeed, if two subgroups $\Gamma', \Gamma'' \subset \Gamma$ give rise to
isometric covers, then $\Gamma'=\gamma^{-1} \Gamma'' \gamma$ for some
$\gamma\in\Isom{\Hyp^n}$ and thus $\gamma$ is in the commensurator of
$\Gamma.$ Moreover, since $\Gamma$ is non-arithmetic, by Margulis's
Theorem(\cite{Ma} Theorem 1, pg. 2),
$[\Comm(\Gamma):\Gamma]=k<\infty$.  Thus a subgroup of $\Gamma$ of
index $r$ is conjugate to at most $rk$ subgroups of $\Gamma$, and thus
we have constructed at least $\frac{(\lfloor{r/2}\rfloor+1)!}{rk}\ge C
e^r$ non-isometric $r$-sheeted covers of diameter at most $D \log{r}$.
Equivalently, there are at least $C e^{e^{d/D}}$ covers of diameter at
most $d$, establishing the lower bound.

\remark This lower bound can also be established by random methods.  Finite index subgroups of a free group $F_k$
correspond to finite covers of a bouquet of $k$ circles, that is, a
$2k$-regular graph.  If $k\ge 5$, then a random $2k$-regular graph is
almost surely an expander and thus has diameter
logarithmic in its number of vertices.  A proof of this fact can be found in \cite{LubExp}.

\section{The Upper Bound}
For the upper bound, note first that in dimension $n\ge 4$, the
argument in \cite{BGLM} holds directly.  If $M$ has diameter $d$, then
$\vol{M}$ is less than the volume of a ball $B^{\Hyp^n}_d$ of radius
$d$ in hyperbolic space, and thus, by the theorem in \cite{BGLM}, for
sufficiently large $d$, 
$$\tau_n(d)\le \rho_n(B^{\Hyp^n}_d) \le \rho_n(c_1e^{(n-1) d}) \le e^{c_2 d e^{(n-1) d}}$$
where $c_1$ and $c_2$ are constants depending on $n$.

The argument of \cite{BGLM} proceeds by covering $M$ with small balls
and using this cover to find $\pi_1(M)$.  Since, by Mostow Rigidity,
the fundamental group of $M$ determines $M$ up to isometry, a bound on
the number of groups obtained this way gives a bound on the number of
possible manifolds.  In the case of $n\ge 4$, it suffices to find the
fundamental group of the thick part of $M$, however, this fails in
dimension 3 because the thin part of the manifold contributes to its
fundamental group.  In this case, because our manifolds have diameter
bounded above and are thus compact, we can find a lower bound on the injectivity radius,
which allows us to use similar methods to find a slightly cruder bound.

If $M$ is a closed hyperbolic manifold with injectivity radius
$\epsilon$, for $\epsilon$ sufficiently small, then its diameter must
be at least on the order of $-c\log{\epsilon}$.  $M$ contains a closed
geodesic $g$ of length $2\epsilon$, and for $\epsilon$ sufficiently
small, $g$ is the central geodesic of a component of the thin part of
$M$.  Reznikov\cite{Rez} shows that the distance between $g$ and the boundary of this component is
at least $-c\log{\epsilon}$, for some constant $c$, giving a lower
bound on the diameter.  Conversely, if $M$ has diameter less than $d$,
then its injectivity radius is at least $e^{-d/c}$.

Let $M$ be a hyperbolic 3-manifold of diameter $d$ and injectivity
radius $r$.  By the above, $r\ge e^{-d/c}$.  If $S$ is a maximal set of points in $M$ such
that any two points in $S$ are separated by at least $r/4$, then the
set $C_M$ of open balls of radius $r/2$ centered at the points of $S$
covers $M$ and the number of points in $S$ is at most
$$\frac{\vol{M}}{\vol{B^{\Hyp^3}_{r/4}}}\le c_1\frac{e^{2d}}{(r/4)^3}
\le c_2e^{5d}.$$ 
Moreover, each of these balls is convex, so intersections of the balls
are convex and thus diffeomorphic to $\R^n$.  By Theorem 13.4 of
\cite{BT}, this implies that the fundamental group of a simplicial 
realization of the nerve of the cover is isomorphic to that of $M$, so
we can calculate $\pi_1(M)$ by considering the combinatorics of this cover.

It suffices to consider the 2-skeleton of the nerve of the cover.  The
1-skeleton is a graph whose vertices are the points of $S$ and whose
edges are the pairs of points in $S$ that are within $r$ of one
another.  If $s\in S$ and $N_s$ is the set of neighbors of $s$, note
that since the points of $N_s$ are separated by $r/4$, the $r/8$-balls
around points of $N_s$ are disjoint and contained in a ball of radius
$r+r/8$ around $s$.  Thus, for $r$ sufficiently small, the degree of
the graph is bounded by a constant $k$ not depending on $r$; the
number of such graphs is at most $|S|^{k|S|}\le e^{c_3de^{5d}}$.

We will bound the number of possible 2-skeleta by considering the number of
triangles in such a graph.  Since each vertex is adjacent to at most
$k$ others, it is a part of at most $k^2$ triangles, for a total of at
most $|S|k^2$ triangles.  The number of possible 2-skeleta is then at
most $|S|^{(k|S|)} 2^{|S|k^2}\le e^{c_4de^{5d}}$, and thus the
number of possible manifolds $M$ is at most $e^{c_4de^{5d}}$. \qed

\end{document}